\documentclass[11pt]{amsart}

\usepackage[dvips]{graphicx}

\usepackage{graphicx}
\usepackage{amssymb}
\usepackage{amsmath}
\usepackage{color}
\usepackage{times}

\usepackage[unicode,bookmarks,colorlinks]{hyperref}
\hypersetup{
    linkcolor=brickred,
}

\definecolor{mahogany}{cmyk}{0, 0.77, 0.87, 0}
\definecolor{salmon}{cmyk}{0, 0.53, 0.38, 0}
\definecolor{melon}{cmyk}{0, 0.46, 0.50, 0}
\definecolor{yellowgreen}{cmyk}{0.44, 0, 0.74, 0}
\definecolor{brickred}{cmyk}{0, 0.89, 0.94, 0.28}
\definecolor{OliveGreen}{cmyk}{0.64, 0, 0.95, 0.40}
\definecolor{RawSienna}{cmyk}{0, 0.72, 1.0, 0.45}
\definecolor{ZurichRed}{rgb}{1, 0, 0} 

\usepackage{amsmath,amstext,amssymb,amsopn,amsthm}
\usepackage{amsmath,amssymb,amsthm}
\usepackage[mathscr]{eucal}

\newtheorem*{thm}{Theorem}

\begin{document}

\newtheorem{remark}{Remark}
\newtheorem{proposition}{Proposition}
\newtheorem{theorem}{Theorem}[section]
\newtheorem{example}{Example}

\newtheorem{corollary}{Corollary}

\numberwithin{definition}{section}
\numberwithin{corollary}{section}


\numberwithin{remark}{section}
\numberwithin{example}{section}
\numberwithin{proposition}{section}

\newcommand{\gap}{\lambda_{2,D}^V-\lambda_{1,D}^V}
\newcommand{\gapR}{\lambda_{2,R}-\lambda_{1,R}}
\newcommand{\bD}{\mathrm{I\! D\!}}
\newcommand{\calD}{\mathcal{D}}
\newcommand{\calA}{\mathcal{A}}
\newcommand{\calM}{\mathbb{M}}
\newcommand{\calL}{\mathbb{L}}

\newcommand{\conjugate}[1]{\overline{#1}}
\newcommand{\abs}[1]{\left| #1 \right|}
\newcommand{\cl}[1]{\overline{#1}}
\newcommand{\expr}[1]{\left( #1 \right)}
\newcommand{\set}[1]{\left\{ #1 \right\}}

\newcommand{\calC}{\mathcal{C}}
\newcommand{\calE}{\mathcal{E}}
\newcommand{\calF}{\mathcal{F}}
\newcommand{\Rd}{\mathbb{R}^d}
\newcommand{\BR}{\mathcal{B}(\Rd)}
\newcommand{\R}{\mathbb{R}}
\newcommand{\al}{\alpha}
\newcommand{\RR}[1]{\mathbb{#1}}
\newcommand{\bR}{\mathrm{I\! R\!}}
\newcommand{\ga}{\gamma}
\newcommand{\om}{\omega}
\newcommand{\A}{\mathbb{A}}
\newcommand{\bH}{\mathbb{H}}
\newcommand{\B}{\mathbb{B}}

\newcommand{\Pro}{\mathbb{P}}
\newcommand\F{\mathcal{F}}
\newcommand\e{\varepsilon}
\def\H{\mathcal{H}}
\def\t{\tau}

\newcommand{\bb}[1]{\mathbb{#1}}
\newcommand{\bI}{\bb{I}}
\newcommand{\bN}{\bb{N}}

\newcommand{\uS}{\mathbb{S}}
\newcommand{\M}{{\mathcal{M}}}
\newcommand{\calB}{{\mathcal{B}}}

\newcommand{\W}{{\mathcal{W}}}

\newcommand{\m}{{\mathcal{m}}}

\newcommand {\mac}[1] { \mathbb{#1} }

\newcommand{\bC}{\Bbb C}

\newcommand{\ang}[1]{\left<#1\right>}  

\newcommand{\brak}[1]{\left(#1\right)}    
\newcommand{\crl}[1]{\left\{#1\right\}}   
\newcommand{\edg}[1]{\left[#1\right]}     

\newcommand{\E}[1]{{\rm E}\left[#1\right]}
\newcommand{\var}[1]{{\rm Var}\left(#1\right)}
\newcommand{\cov}[2]{{\rm Cov}\left(#1,#2\right)}
\newcommand{\N}[1]{||#1||}     
\newcommand{\bM}{\mathbb M}


\title[Trace asymptotics for subordinate semigroups]
{Trace asymptotics for subordinate semigroups}

\author{Rodrigo Ba\~nuelos}\thanks{R. Ba\~nuelos is supported in part  by NSF Grant
\#0603701-DMS}
\address{Department of Mathematics, Purdue University, West Lafayette, IN 47907, USA}
\email{banuelos@math.purdue.edu}
\author{Fabrice Baudoin}\thanks{F. Baudoin is supported in part by NSF Grant \#0907326--DMS}
\address{Department of Mathematics, Purdue University, West Lafayette, IN 47907, USA}
\email{fbaudoin@math.purdue.edu}
\begin{abstract}
We address a conjecture of D. Applebaum on small time trace asymptotics for subordinate Brownian motion on compact manifolds.  
\end{abstract}
\maketitle
\maketitle


In \cite{App1}, heat trace asymptotics  are computed for the semigroup of the square root of the Laplacian (the generator of the Cauchy process) on the $n$-dimensional torus, $SU(2)$ and $SO(3)$, and a conjecture is made \cite[pp. 2493-2494]{App1} that such asymptotics should hold for all $\alpha$-stable processes, $0<\alpha<2$, on an arbitrary compact Lie group.  The purpose of this note is to point out that this is indeed the case and that such results hold for a wide class of subordinations of the Laplacian on compact manifolds.  In fact, as we shall see, such asymptotics follow from Weyl's law \eqref{asymcounting}. We first record some definitions and set some notations.

Let $\mathbb{M}$ be a compact Riemannian manifold of 
dimension $n$ and denote its Riemannian measure by $\mu$. 
We denote the Laplace-Beltrami operator on $\mathbb{M}$ by $-\Delta$ and denote its eigenvalues by 
$
0=\lambda_0 < \lambda_1\le \lambda_2  \le \cdots \to \infty.
$
As is well known, $-\Delta$ generates a heat semigroup $\{\mathbf{P}_t\}$ on $\mathbf{L}^2_\mu (\mathbb{M},\mathbb{R})$ which possesses a heat kernel denoted here by $p(t,x,y)$. Thus 
\begin{equation*}\label{heatdelta}
\mathbf{P}_tf(x)=e^{-t\Delta}f(x)=\int_{\calM}p(t, x, y)f(y)dy,
\end{equation*}
for all $f\in \mathbf{L}^2_\mu (\mathbb{M},\mathbb{R})$. 
 Moreover, there exists a complete orthonormal basis $\{\varphi_k\}_{k \in \mathbb{N}}$ of $\mathbf{L}^2_\mu (\mathbb{M},\mathbb{R})$ consisting of eigenfunctions  $\varphi_k$ having the eigenvalue $\lambda_k$ such that, for $t>0$, $x,y \in \mathbb{M}$,
\begin{equation}\label{eigenexp}
p(t,x,y)=\sum_{k=0}^{\infty} e^{-\lambda_k t} \varphi_k (x) \varphi_k (y),
\end{equation}
with convergence absolute and uniform for each $t>0$.  The operator $\mathbf{P}_t$ is a trace-class  on $\mathbf{L}^2(\bM,\mu)$ and for all $t>0$, 
$\mathbf{Tr} (\mathbf{P}_t)=\sum_{k=0}^{\infty} e^{-\lambda_k t}.$  Let us also denote here the spectral counting function (the number of eigenvalues not exceeding the fixed positive number $\lambda$) of $-\Delta$ by 
$
N(\lambda)=\mathbf{Card} \{ \lambda_i, \lambda_i \le \lambda \}.  
$
The asymptotic of $\mathbf{Tr} (\mathbf{P}_t)$, as $t\to 0$,  and of $N(\lambda)$, as $\lambda\to\infty$, have been extensively studied by many.  In particular, denoting the volume $\mu(\calM)$ of $\calM$ by $Vol(\calM)$, we have  
\begin{equation}\label{asymtrace}
\mathbf{Tr} (\mathbf{P}_t)\sim \frac{Vol(\calM)}{(4\pi)^{n/2}}t^{-n/2}, \qquad t\to  0,  
\end{equation}
and 
\begin{equation}\label{asymcounting}
N(\lambda)\sim \frac{Vol({\calM}) }{\Gamma \left( \frac{n}{2}+1 \right)(4\pi)^{n/2} }\lambda^{n/2},\qquad \lambda\to\infty. 
\end{equation}
Since the trace is the Laplace transform of the spectral counting function, \eqref{asymcounting} immediately implies \eqref{asymtrace}.  On the other hand, from the Karamata tauberian theorem, \eqref{asymcounting} follows from \eqref{asymtrace}, see \cite[p.~109]{Sim}.

The canonical Brownian motion on the manifold $\calM$ with generator $\Delta$ will be denoted by $B_t$.  Let $S_t$ be a subordinator on $[0, \infty)$ independent of $B_t$.  That is, $S_t$ is an increasing L\'evy process taking values in $[0, \infty)$ with $S_0=0$ whose law is determined by its Laplace transform which has the form 
\begin{equation}\label{bern1}
E\left(e^{-\lambda S_t}\right)=e^{-t\psi(\lambda)},\quad  \lambda>0.
\end{equation}
 $\psi:(0, \infty)\to (0, \infty)$ is called the Laplace exponent of $S_t$. It is given by the Bernstein function
\begin{equation}\label{bern2}
\psi(\lambda)=b\lambda +\int_0^{\infty} (1-e^{-\lambda s}) \nu(ds), \quad \lambda>0,
\end{equation}
$b\geq 0$ and  with $\nu$ a positive measure on $(0,\infty)$ (the L\'evy measure of the process) with 
$$
\int_0^{\infty}\left(s\wedge 1\right)\nu(ds) < \infty. 
$$
 For each $t\geq 0$, we let $\eta_t$ be the distribution of $S_t$.  This gives a family of measures which form a convolution semigroup and we can write \eqref{bern1}  as the Laplace transform of the measure $\eta_t$.  That is, 
\begin{equation}\label{bern3}
\int_0^{\infty}e^{-\lambda s}\eta_t(ds)=e^{-t\psi(\lambda)},\quad  \lambda>0. 
\end{equation}
We refer the reader to \cite[Chap. III]{Ber1} for these and other properties of subordinators. 

By the Bochner subordination principle the process $\{B_{S_t}\}$ generates a heat semigroup which we denote here by $\mathbf{Q}_t$ with heat kernel given by 
\begin{equation}\label{subor2}
q(t,x,y) =\int_0^{\infty} p (s,x,y)  \eta_t (ds). 
\end{equation}
By \eqref{eigenexp} and \eqref{bern3} we have that 
$$
q(t, x, y)=\sum_{k=0}^{\infty} e^{-\psi(\lambda_k) t} \varphi_k (x) \varphi_k(y). 
$$
This is the heat kernel for the operator $\psi(-\Delta)$ defined via the spectral theorem. 
Since $\mathbb{M}$ is  compact,  
\[
\int_{\mathbb{M}}\int_{\mathbb{M}} q(t,x,y)^2 d \mu(x) d\mu(y) <\infty.
\]
Thus 
$
\mathbf{Q}_t: \mathbf{L}^2_\mu (\mathbb{M},\mathbb{R}) \rightarrow \mathbf{L}^2_\mu (\mathbb{M},\mathbb{R})
$
and it is a Hilbert-Schmidt operator. In particular, it is compact and  
since $\mathbf{Q}_t=\mathbf{Q}_{t/2} \mathbf{Q}_{t/2}$, $\mathbf{Q}_t$ is a product of two Hilbert-Schmidt operators. It is therefore  class trace and 
\begin{eqnarray*}
\mathbf{Tr} ( \mathbf{Q}_t)& =&\int_{\mathbb{M}}\int_{\mathbb{M}} q(t/2,x,y)q(t/2,y,x) d \mu(x) d\mu(y)
   =\int_{\mathbb{M}} q(t,x,x) d\mu(x) \\
  &=&\sum_{k=0}^{\infty} e^{-\Psi(\lambda_k) t} \int_{\mathbb{M}} \phi_k (x)^2 d\mu (x) 
  =\sum_{k=0}^{\infty} e^{-\Psi(\lambda_n) t} .
\end{eqnarray*}

Examples of subordinators are provided in \cite{Ber1} and  \cite{Cin}, among multiple other sources.  In particular, for any $0<\alpha<2$, set $\psi(\lambda)=\lambda^{\alpha/2}$.   We can write
$$
\lambda^{\alpha/2}=\frac{\alpha/2}{\Gamma(1-\alpha/2)}\int_0^{\infty}(1-e^{-\lambda s})s^{-1-\alpha/2} ds.
$$
This gives the $\alpha$-stable subordinators with L\'evy measures 
$\nu(ds)=\frac{\alpha/2}{\Gamma(1-\alpha/2)}s^{-1-\alpha/2} ds$
and the $\alpha$-stable processes on the manifold $\calM$ with generator $(-\Delta)^{\alpha/2}$ is the Brownian motion subordinated by these processes. What is computed in \cite{App1} is the asymptotics for $\mathbf{Tr} ( \mathbf{Q}_t)$, as $t\downarrow 0$,  for the case $\psi(\lambda)=\sqrt{\lambda}$ when the manifold is the $n$-dimensional torus, $SU(2)$ and $SO(3)$,  and conjectured there also that these asymptotic behavior should hold for any $0<\alpha<2$ for these manifolds.   If we denote by $N^{\psi}(\lambda)$ the spectral counting function of the operator $\mathbf{Q}_t$ associated with $\psi(-\Delta)$ as constructed above, we see that  if $\psi$ is increasing with a well defined increasing inverse denoted by $\phi$ both unbounded, we have 
$N^{\psi}(\lambda)=N(\phi(\lambda))$.  Since the small time asymptotic behavior of $\mathbf{Tr} ( \mathbf{Q}_t)$  can be obtained from the large time behavior of $N^{\psi}(\lambda)$, we see that the behavior found in \cite{App1} for the Cauchy process ($\alpha=1$)  in $SU(2)$ and $SO(3)$ can be extended  to a wide class of subordinators.  
Indeed, if $\psi(\lambda)=\lambda^{\alpha/2}$,  
\begin{equation*}\label{main1}
N^\psi (\lambda)\sim \frac{Vol(\calM)}{\Gamma \left( \frac{n}{2}+1 \right) (4\pi)^{n/2}}\,\lambda^{n/\alpha}, \qquad  \lambda\to\infty
\end{equation*}
and 
\begin{equation*}\label{main2}
\mathbf{Tr} (\mathbf{Q}_t)  =  \frac{Vol(\mathbb{M})\Gamma\left(\frac{n}{\alpha}+1 \right) }{\Gamma \left( \frac{n}{2}+1 \right) (4\pi)^{n/2}}\,t^{-n/\alpha}, \qquad  t\to 0. 
\end{equation*}

From \cite{Ber1},  it follows that if $b>0$ or  $\nu(0, \infty)=\infty$ where $b$ and $\nu$ are as in \eqref{bern2}, then $\psi$ is increasing with $\lim_{+\infty} \psi=+\infty$. Thus $\psi$ has an increasing inverse denoted by $\phi$ and this guarantees that  $N^{\psi}(\lambda)=N(\phi(\lambda))$ and 
\begin{equation}\label{main3}
N^\psi (\lambda)\sim \frac{Vol(\calM)}{\Gamma \left( \frac{n}{2}+1 \right) (4\pi)^{n/2}}\,\phi(\lambda)^{n/2}, \qquad  \lambda\to\infty. 
\end{equation}

We  recall next that a function $L:(0, \infty)\to (0, \infty)$  is said to be  regularly varying at infinity if for all $a>0$, $\lim_{\lambda\to\infty} \frac{L(a\lambda)}{L(\lambda)}=a^{r}$, for some $r\geq 0$. ( If $r=0$, the function is said to be slowly varying at infinity. Thus regularly varying for $r\geq 0$ means that $L(x)=x^{r}l(x)$ where $l$ is slowly varying.)   
We also recall that by the tauberian theorem (see \cite{BinGolTeu}) if $\rho$ is a non-decreasing, right-continuous function on $[0, \infty)$ and 
$\mathfrak{L}{(\rho)}(t) = \int_0^{\infty}{{e^{-st}}d\rho(s)}$
denotes its Laplace transform,  then 
\begin{equation*}
\rho(t) \sim \frac{\gamma}{\Gamma(1 + r)}\, t^{r}\,l(t), \quad t \to \infty \quad \text{iff}\quad 
\mathfrak{L}{\rho}(t) \sim \gamma t^{-r}\,l(1/t), \qquad t \to 0, 
\end{equation*}
whenever $l$ is slowly varying at infinity, $\gamma \geq 0$, and $r \ge 0$.  We also observe  (see \cite{BinGolTeu} Proposition 1.3.6 and Theorem 1.5.12)  that with $\psi$ strictly  increasing and regularly varying  with index $r>0$, then its inverse $\varphi$ is regularly varying with index $\frac{1}{r}$.  In addition, if $\psi$ is regularly varying of index $r>0$, then $\varphi^{n/2}$ is regularly varying of index $\frac{n}{2r}$. 
From these observations we immediately see that \eqref{main3}  imply the asymptotic results for $\mathbf{Tr} (\mathbf{Q}_t)$ for a wideer class of subordinations which include, for example, relativistic Brownian motion and the more general relativistic stable processes.  
(The preceding  reasoning  and calculations are the same as those carried out in \cite{Bry1} and \cite{Bry2} for the Dirichlet problem in Euclidian domains of $\bR^d$. ) The following summarizes the above discussions for the more general subordinations.

\begin{thm}\label{main4}  Suppose the Laplace transform of the subordinator $\psi$ is increasing, regularly varying at infinity with $r>0$ and let $\varphi$ be its inverse. Let $\mathbf{Q}_t$ be the heat semigroup for $\psi(-\Delta)$.  Then 
\begin{equation}\label{main5}
\mathbf{Tr} (\mathbf{Q}_t)  \sim  \frac{Vol(\mathbb{M})\Gamma\left(\frac{n}{2r}+1 \right) }{\Gamma \left( \frac{n}{2}+1 \right) (4\pi)^{n/2}}\,\varphi\left(\frac{1}{t}\right)^{n/2}, \qquad t\to 0. 
\end{equation}
If $r=0$, the same holds replacing $\Gamma\left(\frac{n}{2r}+1 \right)$ with 1. 
\end{thm} 

%
%

Observe also that on a compact Lie group, with normalized Haar measure,  the above trace asymptotics \eqref{main5} implies as a byproduct the on diagonal asymptotics of the heat kernel   
\[
q_t(x,x) \sim \frac{\Gamma\left(\frac{n}{2r}+1 \right) }{\Gamma \left( \frac{n}{2}+1 \right) (4\pi)^{n/2}}\,\varphi\left(\frac{1}{t}\right)^{n/2}, \qquad t\to 0.
\]

\subsection*{Acknowledgments}  We gratefully acknowledge the many useful communications with David Applebaum on the topic of this paper.

\end{document}